\documentclass[12pt]{article}

\usepackage{latexsym}
\usepackage{amsmath}
\numberwithin{equation}{section} 
\usepackage{amsfonts}
\usepackage{amssymb}
\usepackage{amsthm}
\usepackage{bm}
\usepackage{cite} 
\usepackage[colorlinks,linkcolor=red,anchorcolor=blue,citecolor=green,citecolor=green,dvipdfm]{hyperref}
\usepackage{graphicx} 
\usepackage{subfigure} 
\usepackage{color}
\usepackage{caption}
\usepackage{slashbox}  
\usepackage{algorithm} 
\usepackage{algorithmic} 
\usepackage{multirow} 
\usepackage{indentfirst}  
\usepackage[margin=0.9in]{geometry}
\usepackage{rotating}

\newtheorem{theorem}{Theorem}[section]

\newtheorem{definition}[theorem]{Definition} 

\newtheorem{lemma}[theorem]{Lemma} 

\newtheorem{remark}[theorem]{Remark} 

\begin{document}

\title{Low-rank matrix recovery via regularized nuclear norm minimization\footnote{Email addresses: wdwang@swu.edu.cn (Wendong Wang), zhangf@email.swu.edu.cn (Feng Zhang), wjj@swu.edu.cn (Jianjun Wang). The corresponding author is Jianjun Wang.}}

\author{
Wendong Wang$^{1}$, Feng Zhang$^{2}$, and Jianjun Wang$^{2}$
\\
{\small 1. College of Artificial Intelligence, Southwest University, Chongqing, 400715, China}\\
{\small 2. School of Mathematics and Statistics, Southwest University, Chongqing, 400715, China}\\
}

\date{}
\maketitle

\subparagraph{Abstract.}
In this paper, we theoretically investigate the low-rank matrix recovery problem in the context of the unconstrained regularized nuclear norm minimization (RNNM) framework. Our theoretical findings show that, the RNNM method is able to provide a robust recovery of any matrix $X$ (not necessary to be exactly low-rank) from its few noisy measurements $\bm{b}=\mathcal{A}(X)+\bm{n}$ with a bounded constraint $\|\bm{n}\|_{2}\leq\epsilon$, provided that the $tk$-order restricted isometry constant (RIC) of $\mathcal{A}$ satisfies a certain constraint related to $t>0$. Specifically, the obtained recovery condition in the case of $t>4/3$ is found to be same with the sharp condition established previously by Cai and Zhang (2014) to guarantee the exact recovery of any rank-$k$ matrix via the constrained nuclear norm minimization method. More importantly, to the best of our knowledge, we are the first to establish the $tk$-order RIC based coefficient estimate of the robust null space property in the case of $0<t\leq1$.


\section{Introduction}\label{Section1}
Over the past decade, low-rank matrix recovery (LRMR) problem has attracted considerable interest of researchers in many fields, including computer vision \cite{EJ-Candes-RPCA-2011}, recommender systems \cite{R-Mazumder-S-JMLR-2010}, and machine learning \cite{A-Argyriou-C-ML-2008}, to name a few. Mathematically, this problem aims to recover an unknown low-rank matrix $X\in\mathbb{R}^{n_{1}\times n_{2}}$ from
\begin{align*}
\bm{b}=\mathcal{A}(X)+\bm{n},
\end{align*}
where $\bm{b}\in\mathbb{R}^{m}(m\ll n_{1}n_{2})$ is an observed vector, $\bm{n}\in\mathbb{R}^{m}$ is the unknown noise, and $\mathcal{A}: \mathbb{R}^{n_{1}\times n_{2}}\rightarrow \mathbb{R}^{m}$ is a known linear measurement map defined as
\begin{align}\label{Linear-map}
\mathcal{A}(X)=[\text{tr}(X^{T}A^{(1)}), \text{tr}(X^{T}A^{(2)}), \cdots, \text{tr}(X^{T}A^{(m)})]^{T}.
\end{align}
Here, $A^{(i)}$ for $i=1,2,\cdots,m$ is denoted as a matrix with size $n_{1}\times n_{2}$, and $\text{tr}(\cdot)$ is the trace function.

A popular approach for the LRMR problem is to solve a convex nuclear norm minimization (NNM) model
\begin{align}\label{NNM-Model}
\min_{X\in\mathbb{R}^{n_{1}\times n_{2}}}~\|X\|_{*},~~s.t.~~\|\bm{b}-\mathcal{A}(X)\|_{2}\leq\epsilon.
\end{align}
So far, much work has been done to explore the theoretical performance of \eqref{NNM-Model} in exact/robust recovery of any matrix that is not necessary to be exactly
low-rank, see, e.g., \cite{B-Recht-SIAM-Review-2010,EJ-Candes-TIT-2011,foucart2013,MJ-Lai-A-SIAM-2013,delta-sharp-Tony-2013,cai2013compressed,TT-Cai-polytope-TIT-2014,R-Zhang-A-TIT-2018,Zhang-Li-sharp-Lp-InPress}. \textcolor{black}{More specifically, one may seek the sufficient conditions under which the upper-bound estimate of the recovery error will take the form
\begin{align}\label{desired-results-for-constraints}
\|X^{\sharp}-X\|_{F}&\leq C_{1}\frac{\left\|X-X_{[k]}\right\|_{*}}{\sqrt{k}}+ C_{2}\epsilon,
\end{align}
where $X^{\sharp}$ and $X_{[k]}$ are denoted by the optimal solution of \eqref{NNM-Model} and the best rank-$k$ approximate of $X$, respectively, and $C_{1},C_{2}$ are two constants only related to the map $\mathcal{A}$. Note that \eqref{desired-results-for-constraints} also indicates that under these conditions any rank-$k$ matrices, i.e., the matrices whose rank is at most $k$, can be exactly recovered from \eqref{NNM-Model} provided that there is no noise involved, i.e., $\bm{n}$=0 and $\epsilon=0$.} As one of the most powerful and widely used theoretical tools, restricted isometry property (RIP) captures particular attention in establishing these desired conditions and their resulting upper-bound estimates of the recovery error.
\begin{definition}[\cite{EJ-Candes-TIT-2011}]\label{RIP-Definition}
A linear map $\mathcal{A}$ given in \eqref{Linear-map} is said to satisfy the RIP with restricted isometry constant (RIC) of order $k$, denoted by $\delta_{k}$\footnote{When $k$ is not an integer, we define $\delta_{k}$ as $\delta_{\lceil k\rceil}$.}, if $\delta_{k}$ is the smallest value $\delta\in(0,1)$ such
\begin{align*}
(1-\delta)\|X\|_{F}^{2}\leq\|\mathcal{A}(X)\|_{2}^{2}\leq(1+\delta)\|X\|_{F}^{2}
\end{align*}
holds for every rank-$k$ matrix $X\in\mathbb{R}^{n_{1}\times n_{2}}$.
\end{definition}
Some representative conditions include $\delta_{4k}<0.558$ and $\delta_{3k}<0.4721$ in \cite{K-Mohan-NewRIP-ISIT-2010}, $\delta_{2k}<0.4931$ and $\delta_{k}<0.309$ in \cite{wang2013bounds}, and $\delta_{2k}<1/2$ and $\delta_{k}<1/3$ in \cite{delta-sharp-Tony-2013}. In particular, the sharp conditions for exactly rank-$k$ matrix recovery, which takes the form of $\delta_{tk}<\delta^{*}$, have been completely given by Cai and Zhang in \cite{TT-Cai-polytope-TIT-2014} and Zhang and Li in \cite{R-Zhang-A-TIT-2018} for the cases of $0<t<4/3$ and $0<t\leq 4/3$, respectively. To be specific, we can write these sharp conditions into a compact form as below,
\begin{align}\label{sharp-deltaTK}
\delta_{tk}<\left\{
\begin{aligned}
&\frac{t}{4-t}, ~~~~~0<t\leq\frac{4}{3},\\
&\sqrt{\frac{t-1}{t}},  ~~\frac{4}{3}<t<1.
\end{aligned}
\right.
\end{align}
In fact, under the condition \eqref{sharp-deltaTK} any (nearly) low-rank matrix can still be robustly recovered from \eqref{NNM-Model} in the presence of noise, and more details can be found within \cite{TT-Cai-polytope-TIT-2014,R-Zhang-A-TIT-2018}.

\textcolor{black}{Generally, when confronted with the relatively small problems where a high degree of numerical precision is required, one can easily formulate \eqref{NNM-Model} as a semidefinite program (SDP), see, e,.g., \cite{B-Recht-SIAM-Review-2010,EJ-Candes-TIT-2011}, and thus numerically solve it by any of the standard SDP solvers. However, when the scale of the input data is relatively large, it is often not convenient (sometimes maybe impossible) to solve \eqref{NNM-Model} by any standard SDP solvers. Moreover, it is also difficult to estimate a proper parameter value of $\epsilon$ in \eqref{NNM-Model} to well accommodate the unknown noise. Instead of solving \eqref{NNM-Model} directly, many algorithms, see, e.g., \cite{K-Toh-A-PJO-2010,D-Goldfarb-C-FCM-20011,boyd2011distributed,MJ-Lai-IRLS-SIAM-2013}, were proposed to solve the following unconstrained regularized NNM (RNNM) model
\begin{align}\label{NN-Minimization-Unconstrained}
\min_{X\in\mathbb{R}^{n_{1}\times n_{2}}}~\|X\|_{*}+\frac{1}{2\lambda}\|\bm{b}-\mathcal{A}(X)\|_{2}^{2},
\end{align}
where $\lambda>0$ is a trade-off parameter. Compared with the constrained optimization problem \eqref{NNM-Model}, the unconstrained optimization problem \eqref{NN-Minimization-Unconstrained} can well balance the low-rankness of the desired output matrix and the resultant recovery error with properly chosen values of parameter $\lambda$. It has been proved in the practical application that \eqref{NN-Minimization-Unconstrained} is much more suitable for noisy measurements and approximately low-rank matrix recovery \cite{MJ-Lai-IRLS-SIAM-2013}. Nevertheless, one would hope that a result similar to \eqref{desired-results-for-constraints} can be proved for \eqref{NN-Minimization-Unconstrained} as well.} To the best of our knowledge, Cand\`{e}s and Plan \cite{EJ-Candes-TIT-2011} gave the first RIP-based performance guarantee for \eqref{NN-Minimization-Unconstrained}, and their results show that, when the noise $\bm{n}$ obeys $\|\mathcal{A}^{*}(\bm{n})\|\triangleq\|\sum_{i=1}^{m}\bm{n}_{i}\cdot A^{(i)}\|\leq\lambda/2$, and the map $\mathcal{A}$ satisfies $\delta_{4k}<(3\sqrt{2}-1)/17$, the robust recovery of any rank-$k$ matrices can be guaranteed through \eqref{NN-Minimization-Unconstrained}. However, after their initial work, the theoretical investigation of \eqref{NN-Minimization-Unconstrained} is rarely reported. Note that their noise setting is based on the Dantzig selector rather than the often used $\ell_{2}$-norm setting (i.e., $\|\bm{n}\|_{2}\leq\epsilon$), and the obtained sufficient condition still has room to improve.

In this paper, by means of the powerful RIP tool, we theoretically investigate the performance guarantees of the unconstrained RNNM model \eqref{NN-Minimization-Unconstrained} when the noise $\bm{n}$ obeys $\|\bm{n}\|_{2}\leq\epsilon$. In summary, our contributions are two-fold. First, we show that if $\mathcal{A}$ obeys $\delta_{tk}<\sqrt{(t-1)/t}$ for certain $t>1$, then the unconstrained RNNM model \eqref{NN-Minimization-Unconstrained} will be able to provide a robust matrix recovery performance. The obtained sufficient condition is in line with the sharp recovery condition \eqref{sharp-deltaTK} in the case of $t>4/3$ for the constrained problem \eqref{NNM-Model}. Second, by establishing the $tk$-order RIC based coefficient estimate of the robust null space property (RNSP) in the case of $0<t\leq1$, we develop another $tk$-order RIC based sufficient condition for \eqref{NN-Minimization-Unconstrained}, and also obtain some new upper-bound estimates of recovery error.

The remainder of the paper is organized as follows. Section \ref{SECTION-2} introduces some necessary notations and lemmas. Section \ref{SECTION-3} presents a performance guarantee of the RNNM model \eqref{NN-Minimization-Unconstrained} by means of the $tk$-order RIC with $t>1$. In Section \ref{SECTION-4}, we first establish a $tk$-order RIC based coefficient estimate of the RNSP with $0<t\leq1$, and then obtain another parallel performance guarantee result for \eqref{NN-Minimization-Unconstrained}. Finally, conclusion and future work are given in Section \ref{SECTION-5}.

\section{Notations and preliminaries}\label{SECTION-2}
\subsection{Notations}
Without loss of generality we assume that $n_{1}\leq n_{2}$. For any positive integer $k$, we denote $[k]=\{1,2,\cdots, k\}$, and for any $\Omega\subset[n_{1}]$, we denote $\Omega^{c}=[n_{1}]\setminus \Omega$. We denote the singular value decomposition (SVD) of $H\in\mathbb{R}^{n_{1}\times n_{2}}$ as
\begin{align*}
H=\sum_{i=1}^{n_{1}}\sigma_{i}(H)\bm{a}_{H}^{(i)}\left(\bm{c}_{H}^{(i)}\right)^{T},
\end{align*}
where $\sigma_{i}(H)$ is the $i$th largest singular value of $H$, and $\bm{a}_{H}^{(i)}$ and $\bm{c}_{H}^{(i)}$ are the left and right singular value vectors of $H$, respectively. If there is no confusion caused we will write $\sigma_{i}(H)$, $\bm{a}_{H}^{(i)}$ and $\bm{c}_{H}^{(i)}$ as $\sigma_{i}$, $\bm{a}^{(i)}$ and $\bm{c}^{(i)}$ for simplicity, respectively. For convenience, we denote $H^{(i)}=\sigma_{i}\bm{a}^{(i)}\left(\bm{c}^{(i)}\right)^{T}$, $H_{\Omega}=\sum_{i\in \Omega}H^{(i)}$, and also denote by $\sigma_{\Omega}$ the vector whose element is equal to $\sigma_{i}$ for $i\in\Omega$ and 0 otherwise. Then clearly $H_{[k]}=\sum_{i=1}^{k}H^{(i)}$ and $\|\sigma_{\Omega}\|_{1}=\|H_{\Omega}\|_{*}$. In the end, for any given positive number $\alpha$, we denote $T(\alpha,k)\subset \mathbb{R}^{n_{1}}$ as
\begin{align*}
 T(\alpha, k)&=\{\bm{x}\in\mathbb{R}^{n_{1}}: \|\bm{x}\|_{\infty}\leq\alpha, \|\bm{x}\|_{1}\leq k\alpha\},
\end{align*}
and for any $\bm{y}\in\mathbb{R}^{n_{1}}$, we further denote $U(\alpha, k, \bm{y})\subset\mathbb{R}^{n_{1}}$ as
\begin{align*}
 U(\alpha, k, \bm{y})&=\{\bm{x}\in\mathbb{R}^{n_{1}}: \text{supp}(\bm{x})\subseteq\text{supp}(\bm{y}), \|\bm{x}\|_{0}\leq k, \|\bm{x}\|_{1}=\|\bm{y}\|_{1}, \|\bm{x}\|_{\infty}\leq\alpha\},
\end{align*}
where $\|\bm{x}\|_{0}$ is denoted as the number of the nonzero elements in $\bm{x}$.

\subsection{Three key lemmas}
Before presenting our main results, we need some auxiliary lemmas. We start with introducing the first one, which provides a powerful tool to represent a non-sparse vector by the sparse ones. This lemma was first established by Cai and Zhang in \cite{TT-Cai-polytope-TIT-2014}, and later was extended by Zhang and Li in \cite{Zhang-Li-sharp-Lp-InPress}.
\begin{lemma}\label{Sparse-Representation}
Suppose that $\alpha$ is a positive number and $k$ is a positive integer with $k<n_{1}$. Then $\bm{v}\in\mathbb{R}^{n_{1}}$ obeys $\bm{v}\in T(\alpha, k)$ if and only if $\bm{v}$ is in the convex hull of $U(\alpha, k, \bm{v})$. In particular, any $\bm{v}\in T(\alpha, k)$ can be expressed as
\begin{align*}
\bm{v}=\sum_{l}\gamma_{l}\bm{z}^{(l)}
\end{align*}
where $\bm{z}^{(l)}\in U(\alpha, k, \bm{v})$, $0\leq\gamma_{l}\leq1$ and $\sum_{l}\gamma_{l}=1$. Moreover,
\begin{align*}
\sum_{l}\gamma_{l}\|\bm{z}^{(l)}\|_{2}^{2}\leq k\alpha^{2}.
\end{align*}
\end{lemma}
We also need the following Lemma \ref{Lemma-2}, which provides a family of RIC-based conditions under which the RNSP can be guaranteed. More importantly, under these conditions, we will show in Theorem \ref{Theorem-1} that the RNNM model \eqref{NN-Minimization-Unconstrained} is able to robustly recover any matrix that is not necessary to be exactly low-rank.
\begin{lemma}\label{Lemma-2}
For any fixed $t>1$ and any positive integer $k<n_{1}$ with $tk<n_{1}$, if the map $\mathcal{A}$ obeys the RIP of order $tk$ with
\begin{align}\label{Th1_condition}
\delta_{tk}<\sqrt{\frac{t-1}{t}}
\end{align}
then $\mathcal{A}$ have the RNSP with $\beta_{1}>0$ and $0<\beta_{2}<1$. Specifically, for any matrix $H\in\mathbb{R}^{n_{1}\times n_{2}}$ and $\Omega\subset[n_{1}]$ with $|\Omega|=k$, it holds that
\begin{align}\label{Lemma2-results}
\|H_{\Omega}\|_{F}\leq \beta_{1}\|\mathcal{A}(H)\|_{2}+\beta_{2}\frac{\|H_{\Omega^{c}}\|_{*}}{\sqrt{k}},
\end{align}
where
\begin{align*}
\beta_{1}=\frac{2}{(1-\delta_{tk})\sqrt{1+\delta_{tk}}}, \text{~and~} \beta_{2}=\frac{\delta_{tk}}{\sqrt{(1-(\delta_{tk})^{2})(t-1)}}.
\end{align*}
\end{lemma}

\textcolor{black}{
The RNSP, including the classical NSP as a special case, has been demonstrated to be a powerful theoretical tool in providing the robust recovery guarantees of sparse signals or low-rank matrices via some certain constrained optimization problems, see, e.g., \cite{foucart2013,foucart2016flavors,gao2017stability,gao2017recovery}. However, so far it is still an open problem to verify whether a given matrix/map obeys the NSP or not, and also to determine the values of two coefficients (i.e., $\beta_{1}$ and $\beta_{2}$) in RNSP. We note that there exist few researchers who focused on the characterization of these two coefficients with some other theoretical tools (such as the RIC and coherence) that are relatively easy to be checked. To the best of our knowledge, the first $2k$-order RIC based coefficient estimate of RNSP was obtained independently by Shen, et al. in \cite[Lemma 1]{shen2015stable} and Foucart in \cite[Theorem 5]{foucart2016flavors} to fit the sparse recovery scenarios. Later, by using the $tk$-order RIC tool with $t>1$, Ge, et al. in \cite[Lemma 2]{GeHM-UnconstrainedCS-Submitted-2018} extended their results to a more general case. Recently, the coherence-based coefficient estimate of RNSP was obtained by Wang, et al. in \cite[Lemma 3]{wang2019coherence} to deal with the robust signal recovery from the basis pursuit de-noising \cite{Chen-BPDN-SIAM-2001}. In fact, our Lemma \ref{Lemma-2} can be viewed as an extension of \cite[Lemma 2]{GeHM-UnconstrainedCS-Submitted-2018} established for the measurement matrix to that for the measurement map.
}

\begin{proof}[Proof of Lemma \ref{Lemma-2}]
The proof mainly follows from \cite{GeHM-UnconstrainedCS-Submitted-2018}. When $tk$ is not an integer, let $t'=\lceil tk\rceil/k$, then $t'>t$ and $t'k$ is an integer. In view of this, we here only need to prove Lemma \ref{Lemma-2} when $tk$ is a positive integer for a given $t>1$. Let's denote the SVD of $H$ as $H=\sum_{i=1}^{n_{1}}\sigma_{i}\bm{a}^{(i)}\left(\bm{c}^{(i)}\right)^{T}$, and
\begin{align*}
\Lambda_{1}=\left\{i\in\Omega^{c}: \sigma_{i}>\frac{\|H_{\Omega^{c}}\|_{*}}{(t-1)k}\right\},\Lambda_{2}=\left\{i\in \Omega^{c}: \sigma_{i}\leq\frac{\|H_{\Omega^{c}}\|_{*}}{(t-1)k}\right\}.
\end{align*}
Then clearly $\Lambda_{1}\cup \Lambda_{2}=\Omega^{c}$ and $\Lambda_{1}\cap \Lambda_{2}=\emptyset$. In what follows, we start with proving that
\begin{align}\label{Lemma2-proof-inequality1}
\|H_{\Omega\cup \Lambda_{1}}\|_{F}\leq \beta_{1}\|\mathcal{A}(H)\|_{2}+\frac{\beta_{2}}{\sqrt{k}}\|H_{\Omega^{c}}\|_{*}
\end{align}
To do so, we first show that $|\Lambda_{1}|<(t-1)k$. In fact it holds naturally if $\Lambda_{1}=\emptyset$. When $\Lambda_{1}\neq\emptyset$, we know that
\begin{align*}
\|\sigma_{\Lambda_{1}}\|_{1}=\|H_{\Lambda_{1}}\|_{*}>|\Lambda_{1}|\frac{\|H_{\Omega^{c}}\|_{*}}{(t-1)k}\geq\frac{|\Lambda_{1}|}{(t-1)k}\|H_{\Lambda_{1}}\|_{*}
=\frac{|\Lambda_{1}|}{(t-1)k}\|\sigma_{\Lambda_{1}}\|_{1}.
\end{align*}
This also yields the desired result. On the other hand, we can easily induce from the definition of $\Lambda_{1}$ and $\Lambda_{2}$ that
\begin{align*}
\|\sigma_{\Lambda_{2}}\|_{1}&=\|H_{\Omega^{c}}\|_{*}-\|H_{\Lambda_{1}}\|_{*}\leq ((t-1)k-|\Lambda_{1}|)\frac{\|H_{\Omega^{c}}\|_{*}}{(t-1)k},\\
\|\sigma_{\Lambda_{2}}\|_{\infty}&=\max_{i\in \Lambda_{2}}\sigma_{i}\leq\frac{\|H_{\Omega^{c}}\|_{*}}{(t-1)k},
\end{align*}
which, together with Lemma \ref{Sparse-Representation}, indicates that we can express $\sigma_{\Lambda_{2}}$ as
\begin{align*}
\sigma_{\Lambda_{2}}=\sum_{l}\gamma_{l}\bm{z}^{(l)},
\end{align*}
with $\bm{z}^{(l)}$ satisfying
\begin{align}\label{Z-111}
\sum_{l}\gamma_{l}\|\bm{z}^{(l)}\|_{2}^{2}\leq ((t-1)k-|\Lambda_{1}|)\frac{\|H_{\Omega^{c}}\|_{*}^{2}}{(t-1)^{2}k^{2}}\leq \frac{\|H_{\Omega^{c}}\|_{*}^{2}}{(t-1)k}.
\end{align}
By further defining
\begin{align*}
B^{(l)}=(1+\delta_{tk})H_{\Omega\cup \Lambda_{1}}+\delta_{tk}Z^{(l)},~~D^{(l)}=(1-\delta_{tk})H_{\Omega\cup \Lambda_{1}}-\delta_{tk}Z^{(l)},
\end{align*}
where $Z^{(l)}=\sum_{i=1}^{n_{1}}\left(\bm{z}^{(l)}\right)_{i}\bm{a}^{(i)}\left(\bm{c}^{(i)}\right)^{T}$, we can easily induce that both $B^{(l)}$ and $D^{(l)}$ are all rank-$tk$, and $H_{\Lambda_{2}}=\sum_{l}\gamma_{l}Z^{(l)}$. Next, we consider estimating the upper and lower bounds of
\begin{align*}
\rho\triangleq\sum_{l}\gamma_{l}\bigg(\|\mathcal{A}(B^{(l)})\|_{2}^{2}-\|\mathcal{A}(D^{(l)})\|_{2}^{2}\bigg).
\end{align*}
As to the upper bound of $\rho$, we have
\begin{align}\label{Lemma2-proof-inequality3}
\rho&=4\delta_{tk}\langle\mathcal{A}(H_{\Omega\cup \Lambda_{1}}), \mathcal{A}(H_{\Omega\cup \Lambda_{1}}+\sum_{l}\gamma_{i}Z^{(l)})\rangle\nonumber\\
&=4\delta_{tk}\langle\mathcal{A}(H_{\Omega\cup \Lambda_{1}}), \mathcal{A}(H)\rangle\leq4\delta_{tk}\|\mathcal{A}(H_{\Omega\cup \Lambda_{1}})\|_{2}\|\mathcal{A}(H)\|_{2}\nonumber\\
&\leq4\delta_{tk}\sqrt{1+\delta_{tk}}\|H_{\Omega\cup \Lambda_{1}}\|_{F}\|\mathcal{A}(H)\|_{2},
\end{align}
where we have applied the $tk$-order RIP in the last inequality. As to the lower bound of $\rho$, by applying the $tk$-order RIP on $\rho$, we get
\begin{align}\label{Lemma2-proof-inequality4}
\rho&\geq\sum_{l}\gamma_{l}\left((1-\delta_{tk})\|B^{(l)}\|_{F}^{2}-(1+\delta_{tk})\|D^{(l)}\|_{F}^{2}\right)\nonumber\\
&=2\delta_{tk}(1-(\delta_{tk})^2)\|\sigma_{\Omega\cup \Lambda_{1}}\|_{2}^{2}-2(\delta_{tk})^{3}\sum_{l}\gamma_{l}\|\bm{z}^{(l)}\|_{2}^{2}\nonumber\\
&\geq2\delta_{tk}(1-(\delta_{tk})^2)\|H_{\Omega\cup \Lambda_{1}}\|_{F}^{2}-\frac{2(\delta_{tk})^{3}}{(t-1)k}\|H_{\Omega^{c}}\|_{*}^{2},
\end{align}
where we have used $\langle\sigma_{\Omega\cup \Lambda_{1}},\bm{z}^{(l)}\rangle=0$ in the first equality and \eqref{Z-111} in the last inequality. Therefore, combing \eqref{Lemma2-proof-inequality3} and \eqref{Lemma2-proof-inequality4} gives
\begin{align*}
(1-(\delta_{tk})^2)\|H_{\Omega\cup \Lambda_{1}}\|_{F}^{2}-2\sqrt{1+\delta_{tk}}\|\mathcal{A}(H)\|_{2}\|H_{\Omega\cup \Lambda_{1}}\|_{F}-\frac{(\delta_{tk})^{2}}{(t-1)k}\|H_{\Omega^{c}}\|_{*}^{2}\leq0.
\end{align*}
Therefore,
\begin{align*}
\|H_{E\cup E_{1}}\|_{F}\leq&\frac{2\sqrt{1+\delta_{tk}}\|\mathcal{A}(H)\|_{2}}{2(1-(\delta_{tk})^2)}\\
&+\frac{\sqrt{(2\sqrt{1+\delta_{tk}}\|\mathcal{A}(H)\|_{2})^{2}+4(1-(\delta_{tk})^{2})\frac{(\delta_{tk})^{2}}{(t-1)k}\|H_{E^{c}}\|_{*}^{2}}}{2(1-(\delta_{tk})^2)}\\
\leq&\frac{2(1-\delta_{tk})^{-1}}{\sqrt{1+\delta_{tk}}}\|\mathcal{A}(H)\|_{2}+\frac{\delta_{tk}}{\sqrt{(1-(\delta_{tk})^{2})(t-1)}}
\frac{\|H_{\Omega^{c}}\|_{*}}{\sqrt{k}},
\end{align*}
where we have used $\sqrt{x^{2}+y^{2}}\leq|x|+|y|$ for any $x,y\in\mathbb{R}$ in the last inequality. This, together with $\|H_{\Omega}\|_{F}\leq\|H_{\Omega\cup \Lambda_{1}}\|_{F}$, directly yields \eqref{Lemma2-results}. The obtained condition \eqref{Th1_condition} follows trivially from \eqref{Lemma2-results} by enforcing $\beta_{2}=\delta_{tk}/\sqrt{(1-(\delta_{tk})^{2})(t-1)}<1$.
\end{proof}

In the end, we introduce the last lemma (i.e., Lemma \ref{Lemma-2}), which characterizes the relationship between the original solution $X$ and the optimal solution $X^{\sharp}$ of \eqref{NN-Minimization-Unconstrained}.
\begin{lemma}\label{Lemma-3}
Assume that $X^{\sharp}$ is the solution of \eqref{NN-Minimization-Unconstrained} and $H=X^{\sharp}-X$. If the noisy measurements $\bm{b}=\mathcal{A}(X)+\bm{n}$ are observed with the noise level $\|\bm{n}\|_{2}\leq\epsilon$, then for any subset $\Omega\subset[n_{1}]$ with $|\Omega|=k$, we have
\begin{align}\label{Lemma3-results1}
\begin{split}
\|\mathcal{A}(H)\|_{2}^{2}-2\epsilon\|\mathcal{A}(H)\|_{2}\leq&2\lambda(\|H_{\Omega}\|_{*}-\|H_{\Omega^{c}}\|_{*}+2\|X_{\Omega^{c}}\|_{*})
\end{split}
\end{align}
and
\begin{align}\label{Lemma3-results2}
\|H_{\Omega^{c}}\|_{*}\leq\|H_{\Omega}\|_{*}+2\|X_{\Omega^{c}}\|_{*}+\frac{\epsilon}{\lambda}\|\mathcal{A}(H)\|_{2}.
\end{align}
\end{lemma}
\begin{proof}[Proof of Lemma \ref{Lemma-3}]
Since $X^{\sharp}$ is the optimal solution of \eqref{NN-Minimization-Unconstrained}, we have
\begin{align*}
\|X^{\sharp}\|_{*}+\frac{1}{2\lambda}\|\bm{b}-\mathcal{A}(X^{\sharp})\|_{2}^{2}\leq\|X\|_{*}+\frac{1}{2\lambda}\|\bm{b}-\mathcal{A}(X)\|_{2}^{2},
\end{align*}
which is equivalent to
\begin{align}\label{Estimate-LR}
\|\mathcal{A}(H)\|_{2}^{2}-2\langle \bm{n}, \mathcal{A}(H)\rangle\leq 2\lambda(\|X\|_{*}-\|X^{\sharp}\|_{*}).
\end{align}
As to the left-hand side (LHS) of \eqref{Estimate-LR}, we have
\begin{align}\label{Estimate-L}
\|\mathcal{A}(H)\|_{2}^{2}-2\langle \bm{n}, \mathcal{A}(H)\rangle\geq\|\mathcal{A}(H)\|_{2}^{2}-2\epsilon\|\mathcal{A}(H)\|_{2}.
\end{align}
As to the right-hand side (RHS) of \eqref{Estimate-LR}, we know
\begin{align}\label{Estimate-R}
\|X^{\sharp}\|_{*}-\|X\|_{*}=&\sum_{i=1}^{n_{1}}\sigma_{i}(X+H)-(\|X_{\Omega}\|_{*}+\|X_{\Omega^{c}}\|_{*})\nonumber\\
\geq&\sum_{i=1}^{n_{1}}|\sigma_{i}(X)-\sigma_{i}(-H)|-(\|X_{\Omega}\|_{*}+\|X_{\Omega^{c}}\|_{*})\nonumber\\
\geq&\sum_{i\in \Omega}(\sigma_{i}(X)-\sigma_{i}(H))+\sum_{i\in \Omega^{c}}(\sigma_{i}(H)-\sigma_{i}(X))-(\|X_{\Omega}\|_{*}+\|X_{\Omega^{c}}\|_{*})\nonumber\\
=&-\|H_{\Omega}\|_{*}+\|H_{E^{c}}\|_{*}-2\|X_{\Omega^{c}}\|_{*},
\end{align}
where we have used \cite[Theorem 1]{MC-Yue-A-ACHA-2016} in the first inequality. Then combing \eqref{Estimate-LR}, \eqref{Estimate-L}, and \eqref{Estimate-R} leads to the desired result \eqref{Lemma3-results1}, and \eqref{Lemma3-results2} follows trivially from \eqref{Lemma3-results1}.
\end{proof}

\section{\textcolor{black}{Performance guarantee of RNMM model under $tk$-order RIC with $t>1$}}\label{SECTION-3}
With previous preparations in mind, we present our first theoretical result.
\begin{theorem}\label{Theorem-1}
For any observed vector $\bm{b}=\mathcal{A}(X)+\bm{n}$ with a bounded constraint $\|\bm{n}\|_{2}\leq\epsilon$, if the $tk$-order RIC of $\mathcal{A}$ with $t>1$ satisfies condition \eqref{Th1_condition}, then we have
\begin{align}
\|\mathcal{A}(X^{\sharp}-X)\|_{2}&\leq C_{1}(\beta_{1},\beta_{2})\|X-X_{[k]}\|_{*}+ C_{2}(\beta_{1},\beta_{2}), \label{Theorem-1-results1}\\
\|X^{\sharp}-X\|_{F}&\leq C_{3}(\beta_{1},\beta_{2})\|X-X_{[k]}\|_{*}+ C_{4}(\beta_{1},\beta_{2}),\label{Theorem-1-results2}
\end{align}
where $X^{\sharp}$ is the optimal solution of \eqref{NN-Minimization-Unconstrained}, and
\begin{align*}
C_{1}(\beta_{1},\beta_{2})&=\frac{2\lambda}{\sqrt{k}\beta_{1}\lambda+\epsilon},~~C_{2}(\beta_{1},\beta_{2})=2\left(\sqrt{k}\beta_{1}\lambda+\epsilon\right),\\
C_{3}(\beta_{1},\beta_{2})&=\frac{2\sqrt{k}\beta_{1}[3+3\beta_{2}+\left(\beta_{2}\right)^{2}]\lambda + 2[1+4\beta_{2}+2\left(\beta_{2}\right)^{2}]\epsilon }{\sqrt{k}(1-\beta_{2})(\sqrt{k}\beta_{1}\lambda+\epsilon)},\\
C_{4}(\beta_{1},\beta_{2})&=\frac{\sqrt{k}\beta_{1}(5+2\beta_{2})\lambda+
[1+4\beta_{2}+2(\beta_{2})^{2}]\epsilon}{\sqrt{k}(1-\beta_{2})\lambda}\left(\sqrt{k}\beta_{1}\lambda+\epsilon\right).
\end{align*}
\end{theorem}
\textcolor{black}{
\begin{remark}
The condition \eqref{Th1_condition} has been obtained previously by Cai and Zhang in \cite{TT-Cai-polytope-TIT-2014} for exact/robust signal recovery from \eqref{NNM-Model}, and it has been proved to be sharp for the exactly rank-$k$ matrix recovery when $t>4/3$. To the best of our knowledge, we first extend nontrivially this condition from the constrained NNM model \eqref{NNM-Model} to its unconstrained counterpart, i.e., the unconstrained RNNM model \eqref{NN-Minimization-Unconstrained}. On the other hand, note that the obtained coefficients $C_{i}(\beta_{1},\beta_{2})$ (for $i=1,2,3,4$) might seem a bit complicated since they not only involve $\beta_{1},\beta_{2}$, but also involve $k$, $\lambda$ and $\epsilon$. To remedy this problem, we need to do some simplification. We here only take $C_{3}(\beta_{1},\beta_{2})$ and $C_{4}(\beta_{1},\beta_{2})$ for examples. Since
\begin{align*}
C_{3}(\beta_{1},\beta_{2})\leq\frac{2\beta_{1}[3+3\beta_{2}+\left(\beta_{2}\right)^{2}] + 2[1+4\beta_{2}+2\left(\beta_{2}\right)^{2}](\epsilon/\lambda) }{\sqrt{k}\beta_{1}(1-\beta_{2})}
\end{align*}
and
\begin{align*}
C_{4}(\beta_{1},\beta_{2})\leq\frac{\beta_{1}(5+2\beta_{2})+
[1+4\beta_{2}+2(\beta_{2})^{2}](\epsilon/\lambda)}{(1-\beta_{2})}\cdot\sqrt{k}\left[\beta_{1}+(\epsilon/\lambda)\right]\lambda,
\end{align*}
we thus can induce from \eqref{Theorem-1-results2} that
\begin{align}\label{New-FF}
\|X^{\sharp}-X\|_{F}&\leq \widehat{C}_{\lambda/\epsilon}\frac{\|X-X_{[k]}\|_{*}}{\sqrt{k}} + \sqrt{k}\widetilde{C}_{\lambda/\epsilon}\lambda,
\end{align}
where $\widehat{C}_{\lambda/\epsilon}$ and $\widetilde{C}_{\lambda/\epsilon}$ are two constants only relying on the map $\mathcal{A}$ and the value of $\lambda/\epsilon$, and they are given as below.
\begin{align*}
\widehat{C}_{\lambda/\epsilon}=&\frac{2\beta_{1}[3+3\beta_{2}+\left(\beta_{2}\right)^{2}] + 2[1+4\beta_{2}+2\left(\beta_{2}\right)^{2}](\epsilon/\lambda) }{\beta_{1}(1-\beta_{2})},\\
\widetilde{C}_{\lambda/\epsilon}=&\frac{\beta_{1}(5+2\beta_{2})+
[1+4\beta_{2}+2(\beta_{2})^{2}](\epsilon/\lambda)}{(1-\beta_{2})\left[\beta_{1}+(\epsilon/\lambda)\right]^{-1}}.
\end{align*}
Note that the induced upper-bound estimate \eqref{New-FF} also coincide with the ones established in \cite{shen2015stable,GeHM-UnconstrainedCS-Submitted-2018,wang2019coherence,li2019signal} in form.
\end{remark}
}
\textcolor{black}{
\begin{remark}
Theorem \ref{Theorem-1} states that if the measurement map $\mathcal{A}$ obeys a certain $tk$-order RIP condition related to $t>1$, any matrix that is not necessary to be exactly low-rank can be robustly recovered from \eqref{NN-Minimization-Unconstrained} for any fixed parameter $\lambda>0$. According to the obtained results, it is difficult to determine a ``good'' parameter $\lambda$ to yield a ``good'' solution in general case. In fact, so far it is still an open problem to theoretically determine a general parameter $\lambda$ to make sure
that the unconstrained RNNM model \eqref{NN-Minimization-Unconstrained} can perform well. However, if taking a close look at the obtained \eqref{Theorem-1-results2} and \eqref{New-FF}, one will find that the selected parameter $\lambda$ should not be much too large or small. Furthermore, if the desired matrix $X$ is assumed to be exactly rank-$k$, then we can induce from \eqref{Theorem-1-results2} that
\begin{align*}
\|X^{\sharp}-X\|_{F}\leq  C_{4}(\beta_{1},\beta_{2}).
\end{align*}
Obviously, if we desire a optimal solution with recovery error as small as possible from \eqref{NN-Minimization-Unconstrained}, we need to make sure that the value of $C_{4}(\beta_{1},\beta_{2})$ is also as small as possible with respect to the parameter $\lambda$. Considering that
\begin{align*}
C_{4}(\beta_{1},\beta_{2})=&\frac{1}{\sqrt{k}(1-\beta_{2})}\bigg\{k(\beta_{1})^{2}(5+2\beta_{2})\lambda+\frac{1}{\lambda}
\left[1+4\beta_{2}+\left(\beta_{2}\right)^{2}\right]\epsilon^{2}\\
&+2\sqrt{k}\beta_{1}\left[3+3\beta_{2}+\left(\beta_{2}\right)^{2}\right]\epsilon\bigg\}\\
\overset{(a)}{\geq}&\frac{2\beta_{1}}{(1-\beta_{2})}\bigg\{\sqrt{(5+2\beta_{2})\left[1+4\beta_{2}+\left(\beta_{2}\right)^{2}\right]}
+\left[3+3\beta_{2}+\left(\beta_{2}\right)^{2}\right]\bigg\}\epsilon,
\end{align*}
where the equality in (a) holds when $\lambda$ satisfies
\begin{align}\label{lambda-best}
\lambda=\frac{1+4\beta_{2}+\left(\beta_{2}\right)^{2}}{k\left(\beta_{1}\right)^{2}(5+2\beta_{2})}\epsilon^{2},
\end{align}
a ideal selection of parameter $\lambda$ is to set it as in \eqref{lambda-best}, which is related to the two coefficient estimates of RNSP of the map $\mathcal{A}$, noise level $\epsilon$ and also the rank parameter $k$. In realistic situations, such a setting of $\lambda$ is impractical. However, from \eqref{lambda-best} we can capture some information to help set a proper $\lambda$, i.e., the value of $\lambda$ is proportional to that of $\epsilon^{2}$, and inversely proportional to that of $k$.
\end{remark}
}

Now, we present the proof of Theorem \ref{Theorem-1} as follows.
\begin{proof}
We start with proving \eqref{Theorem-1-results1}. Let's define $\widehat{\Omega}=[k]$ and $H=X^{\sharp}-X$. Then by using Lemma \ref{Lemma-2} and Lemma \ref{Lemma-3} with $\Omega=\widehat{\Omega}$, we have
\begin{align}\label{Theorem1-proof-inequlity1}
\|\mathcal{A}(H)\|_{2}^{2}-2\epsilon\|\mathcal{A}(H)\|_{2}\leq&2\lambda(\sqrt{k}\|H_{\widehat{\Omega}}\|_{F}-\|H_{\widehat{\Omega}^{c}}\|_{*}+2\|X_{\widehat{\Omega}^{c}}\|_{*})\nonumber\\
\leq&2\sqrt{k}\lambda(\beta_{1}\|\mathcal{A}(H)\|_{2}+\frac{\beta_{2}}{\sqrt{k}}\|H_{\widehat{\Omega}^{c}}\|_{*})\nonumber\\
&-2\lambda\|H_{\widehat{\Omega}^{c}}\|_{*}+4\lambda\|X_{\widehat{\Omega}^{c}}\|_{*}\nonumber\\
=&2\sqrt{k}\beta_{1}\lambda\|\mathcal{A}(H)\|_{2}-2(1-\beta_{2})\lambda\|H_{\widehat{\Omega}^{c}}\|_{*}+4\lambda\|X_{\widehat{\Omega}^{c}}\|_{*}.
\end{align}
Due to \eqref{Th1_condition}, $\beta_{2}<1$ and thus we induce from \eqref{Theorem1-proof-inequlity1} that
\begin{align*}
\|\mathcal{A}(H)\|_{2}^{2}-2(\sqrt{k}\beta_{1}\lambda+\epsilon)\|\mathcal{A}(H)\|_{2}-4\lambda\|X_{\widehat{\Omega}^{c}}\|_{*}\leq0.
\end{align*}
This directly leads to
\begin{align*}
\|\mathcal{A}(H)\|_{2}\leq&(\sqrt{k}\beta_{1}\lambda+\epsilon)+\sqrt{(\sqrt{k}\beta_{1}\lambda+\epsilon)^{2}+ 4\lambda\|X_{\widehat{\Omega}^{c}}\|_{*}}\\
\leq&(\sqrt{k}\beta_{1}\lambda+\epsilon)+(\sqrt{k}\beta_{1}\lambda+\epsilon)+\frac{2\lambda\|X_{\widehat{\Omega}^{c}}\|_{*}}{(\sqrt{k}\beta_{1}\lambda+\epsilon)}\\
\leq&\frac{2\lambda}{\sqrt{k}\beta_{1}\lambda+\epsilon}\|X_{\widehat{\Omega}^{c}}\|_{*}+2\sqrt{k}\beta_{1}\lambda+2\epsilon,
\end{align*}
which is the desired \eqref{Theorem-1-results1}.

Before proving \eqref{Theorem-1-results2}, let's define $\Omega_{1}=\{k+1,k+2,\cdots, 2k\}$, $\Omega_{2}=\{2k+1,2k+2,\cdots, 3k\}$, $\Omega_{3}=\{3k+1,3k+2,\cdots, 4k\}$, and so on. Thus for $i=2,3,4,\cdots$, we have $\|H_{\Omega_{i}}\|_{F}\leq\|H_{\Omega_{1}}\|_{F}$, and therefore,
\begin{align}\label{previous-results-2015}
\|H_{\widehat{\Omega}^{c}}\|_{F}=&\sqrt{\left\|H_{\Omega_{1}}\right\|_{F}^{2}+\sum_{i\geq2}\left\|H_{\Omega_{i}}\right\|_{F}^{2}}\leq\sqrt{\left\|H_{\Omega_{1}}\right\|_{F}^{2}+\left\|H_{\Omega_{1}}\right\|_{F}\sum_{i\geq2}
\left\|H_{\Omega_{i}}\right\|_{F}}\nonumber\\
\leq&\left\|H_{\Omega_{1}}\right\|_{F}+\frac{1}{2}\sum_{i\geq2}\left\|H_{\Omega_{i}}\right\|_{F}\leq\left\|H_{\Omega_{1}}\right\|_{F}+\frac{1}{2\sqrt{k}}\left\|H_{\widehat{\Omega}^{c}}\right\|_{*},
\end{align}
where the last inequality is due the fact that
\begin{align*}
\sum_{i\geq2}\left\|H_{\Omega_{i}}\right\|_{F}\leq\sqrt{k}\sum_{i\geq2}\max_{j\in\Omega_{i}}\sigma_{j}\leq\frac{1}{\sqrt{k}}\sum_{i\geq2}\|H_{\Omega_{i-1}}\|_{*}=\frac{1}{\sqrt{k}}\|H_{\widehat{\Omega}^{c}}\|_{*}.
\end{align*}
Note that by combining Lemma \ref{Lemma-2} and Lemma \ref{Lemma-3} with $\Omega=\widehat{\Omega}$ again, we can provide two upper-bound estimates of $\|H_{\widehat{\Omega}}\|_{F}$ and $\|H_{\widehat{\Omega}^{c}}\|_{*}$, respectively, which are independent from each other. First, as to that of $\|H_{\widehat{\Omega}}\|_{F}$, we have
\begin{align*}
\|H_{\widehat{\Omega}}\|_{F}\leq&\beta_{1}\|\mathcal{A}(H)\|_{2}+\frac{\beta_{2}}{\sqrt{k}}\left(\|H_{\widehat{\Omega}}\|_{*}+2\|X_{\widehat{\Omega}^{c}}\|_{*}+\frac{\epsilon}{\lambda}\|\mathcal{A}(H)\|_{2}\right)\\
\leq&\beta_{1}\|\mathcal{A}(H)\|_{2}+\beta_{2}\|H_{\widehat{\Omega}}\|_{F}+\frac{\beta_{2}}{\sqrt{k}}\left(2\|X_{\widehat{\Omega}^{c}}\|_{*}+\frac{\epsilon}{\lambda}\|\mathcal{A}(H)\|_{2}\right),
\end{align*}
which is equivalent to
\begin{align}\label{HOmega1}
\|H_{\widehat{\Omega}}\|_{F}\leq&\frac{\beta_{1}}{1-\beta_{2}}\|\mathcal{A}(H)\|_{2}+\frac{\beta_{2}}{\sqrt{k}(1-\beta_{2})}\left(2\|X_{\widehat{\Omega}^{c}}\|_{*}+\frac{\epsilon}{\lambda}\|\mathcal{A}(H)\|_{2}\right)\nonumber\\
=&\frac{\sqrt{k}\beta_{1}\lambda+\beta_{2}\epsilon}{\sqrt{k}(1-\beta_{2})\lambda}\|\mathcal{A}(H)\|_{2}+\frac{2\beta_{2}}{\sqrt{k}(1-\beta_{2})}\|X_{\widehat{\Omega}^{c}}\|_{*}.
\end{align}
Similarly, we can also easily get the upper-bound estimate of $\|H_{\widehat{\Omega}^{c}}\|_{*}$ as below.
\begin{align}\label{HOmegaC1}
\|H_{\widehat{\Omega}^{c}}\|_{*}\leq\frac{\sqrt{k}\beta_{1}\lambda+\epsilon}{(1-\beta_{2})\lambda}\|\mathcal{A}(H)\|_{2}+\frac{2}{1-\beta_{2}}\|X_{\widehat{\Omega}^{c}}\|_{*}.
\end{align}

On the other hand, by combining \eqref{previous-results-2015} and Lemma \ref{Lemma-2} with $\Omega=\Omega_{1}$, we have
\begin{align*}
\|H_{\Omega_{1}}\|_{F}\leq&\beta_{1}\|\mathcal{A}(H)\|_{2}+\beta_{2}\frac{\|H_{\Omega_{1}^{c}}\|_{*}}{\sqrt{k}}\leq \beta_{1}\|\mathcal{A}(H)\|_{2}+\frac{\beta_{2}}{\sqrt{k}}\left(\|H_{\widehat{\Omega}}\|_{*}+\|H_{\widehat{\Omega}^{c}}\|_{*}\right),\\
\leq&\beta_{1}\|\mathcal{A}(H)\|_{2}+\beta_{2}\left(\|H_{\widehat{\Omega}}\|_{F}+\frac{\|H_{\widehat{\Omega}^{c}}\|_{*}}{\sqrt{k}}\right),
\end{align*}
which, together with \eqref{HOmega1} and \eqref{HOmegaC1}, yields
\begin{align}\label{HOmega11}
\|H_{\Omega_{1}}\|_{F}\leq&\beta_{1}\|\mathcal{A}(H)\|_{2}+\beta_{2}\|H_{\widehat{\Omega}}\|_{F}+\frac{\beta_{2}}{\sqrt{k}}\|H_{\widehat{\Omega}^{c}}\|_{*}\nonumber\\
\leq&\beta_{1}\|\mathcal{A}(H)\|_{2}
+\beta_{2}\left(\frac{\sqrt{k}\beta_{1}\lambda+\beta_{2}\epsilon}{\sqrt{k}(1-\beta_{2})\lambda}\|\mathcal{A}(H)\|_{2}+\frac{2\beta_{2}}{\sqrt{k}(1-\beta_{2})}\|X_{\widehat{\Omega}^{c}}\|_{*}\right)\nonumber\\
&+\frac{\beta_{2}}{\sqrt{k}}\left(\frac{\sqrt{k}\beta_{1}\lambda+\epsilon}{(1-\beta_{2})\lambda}\|\mathcal{A}(H)\|_{2}+\frac{2}{1-\beta_{2}}\|X_{\widehat{\Omega}^{c}}\|_{*}\right)\nonumber\\
\leq& \frac{(1+\beta_{2})(\sqrt{k}\beta_{1}\lambda+\beta_{2}\epsilon) }{ \sqrt{k}(1-\beta_{2})\lambda }\|\mathcal{A}(H)\|_{2} + \frac{2\beta_{2}(1+\beta_{2})}{\sqrt{k}(1-\beta_{2})} \|X_{\widehat{\Omega}^{c}}\|_{*}.
\end{align}
Now by combining \eqref{previous-results-2015}, \eqref{HOmega1}, \eqref{HOmegaC1} and \eqref{HOmega11}, we can estimate $\|H\|_{F}$ as follows.
\begin{align*}
\|H\|_{F}\leq&\|H_{\widehat{\Omega}}\|_{F}+\|H_{\widehat{\Omega}^{c}}\|_{F}\leq\|H_{\widehat{\Omega}}\|_{F}
+\left\|H_{\Omega_{1}}\right\|_{F}+\frac{1}{2\sqrt{k}}\left\|H_{\widehat{\Omega}^{c}}\right\|_{*}\\
\leq&\frac{\sqrt{k}\beta_{1}(5+2\beta_{2})\lambda+[1+4\beta_{2}+2(\beta_{2})^{2}]\epsilon}{2\sqrt{k}(1-\beta_{2})\lambda}\|\mathcal{A}(H)\|_{2}
+\frac{1+2\beta_{2}(2+\beta_{2})}{\sqrt{k}(1-\beta_{2})}\|X_{\widehat{\Omega}^{c}}\|_{*},
\end{align*}
which, together with \eqref{Theorem-1-results1}, yields
\begin{align*}
\|H\|_{F}\leq&\frac{\sqrt{k}\beta_{1}(5+2\beta_{2})\lambda+[1+4\beta_{2}+2(\beta_{2})^{2}]\epsilon}{2\sqrt{k}(1-\beta_{2})\lambda}
\left(\frac{2\lambda}{\sqrt{k}\beta_{1}\lambda+\epsilon}\|X_{\widehat{\Omega}^{c}}\|_{*}+2\sqrt{k}\beta_{1}\lambda+2\epsilon\right)\\
&+\frac{1+2\beta_{2}(2+\beta_{2})}{\sqrt{k}(1-\beta_{2})}\|X_{\widehat{\Omega}^{c}}\|_{*}\\
\leq&\frac{2\sqrt{k}\beta_{1}[3+3\beta_{2}+\left(\beta_{2}\right)^{2}]\lambda + 2[1+4\beta_{2}+2\left(\beta_{2}\right)^{2}]\epsilon }{\sqrt{k}(1-\beta_{2})(\sqrt{k}\beta_{1}\lambda+\epsilon)}\|X_{\widehat{\Omega}^{c}}\|_{*}\\
&+\frac{\sqrt{k}\beta_{1}(5+2\beta_{2})\lambda+[1+4\beta_{2}+2(\beta_{2})^{2}]\epsilon}{\sqrt{k}(1-\beta_{2})\lambda}\left(\sqrt{k}\beta_{1}\lambda+\epsilon\right).
\end{align*}
This completes the proof of \eqref{Theorem-1-results2}.
\end{proof}

\section{\textcolor{black}{$tk$-order RIC based coefficient estimate of RNSP with $0<t\leq1$}}\label{SECTION-4}
In the previous section, a family of $tk$-order RIP conditions and their resultant recovery error estimate results are established for the robust matrix recovery from the unconstrained RNNM model \eqref{NN-Minimization-Unconstrained}. As is seen from Theorem \ref{Theorem-1} and its proof, Lemma \ref{Lemma-2}, i.e., the RNSP with $tk$-order RIC based coefficient estimate, plays a vital role in establishing the desired results. Unfortunately, Lemma \ref{Lemma-2}, as well as its resultant Theorem \ref{Theorem-1}, only considers the case of $t>1$. In this section, we will show that under the $tk$-order RIP condition with $0<t\leq1$, \eqref{NN-Minimization-Unconstrained} is still able to provide a robust recovery performance. Before moving on, we have to introduce \cite[Lemma 1]{R-Zhang-A-TIT-2018} since it will be frequently used in the proof of our main results, and one can find it in Lemma \ref{Key-Combination-Results}.
\begin{lemma}\label{Key-Combination-Results}
Let $\bm{w}\in\mathbb{R}^{k}$ be a vector with $\bm{w}=[w_{1},w_{2},\cdots, w_{k}]^{T}$. Choose all subsets $S_{i}\subseteq[k]$ with $|S_{i}|=s<k$, $i\in I$ with $|I|=\binom{k}{s}$, then we have
\begin{align*}
\sum_{i\in I}\sum_{l\in S_{i}}w_{l}=\binom{k-1}{s-1}\sum_{l}w_{l} ~~(s\geq1),
\end{align*}
and
\begin{align*}
\sum_{i\in I}\sum_{\substack{p\neq q\\p,q\in S_{i}}}w_{p}w_{q}=\binom{k-2}{s-2}\sum_{p\neq q}w_{p}w_{q}~~(s\geq2).
\end{align*}
\end{lemma}
Now we are ready to present our second theoretical result, i.e., the $tk$-order RIC based coefficient estimate of RNSP with $0<t\leq1$.
\begin{theorem}\label{key-Theorem-2}
For any fixed $0<t\leq1$ and any positive integer $k<n_{1}$ with $tk<n_{1}$, if the map $\mathcal{A}$ obeys the RIP of order $tk$ with
\begin{align}\label{conditions-step1}
\delta_{tk}\leq\left\{
\begin{aligned}
&\frac{(2-t)t}{(2-t)^{2}+2t\sqrt{(2-t)(1-t)}+(2-t)\theta_{1}}, ~~0<t\leq\frac{2}{3},\\
&\frac{t^{2}}{2(t+1)},  ~~\frac{2}{3}<t\leq1,
\end{aligned}
\right.
\end{align}
where $\theta_{1}=\frac{1}{t}\sqrt{\frac{2-t}{1-t}}$, then $\mathcal{A}$ obeys the RNSP with $\widehat{\beta}_{1}>0$ and $0<\widehat{\beta}_{2}<1$. Specifically, for any matrix $H\in\mathbb{R}^{n_{1}\times n_{2}}$ and any subset $\Omega\subset[n_{1}]$ with $|\Omega|=k$, it holds that
\begin{align}\label{Main-results-2}
\|H_{\Omega}\|_{F}\leq& \widehat{\beta}_{1}\|\mathcal{A}(H)\|_{2}+\widehat{\beta}_{2}\frac{\|H_{\Omega^{c}}\|_{*}}{\sqrt{k}},
\end{align}
where
\begin{align*}
\widehat{\beta}_{1}=\left\{
\begin{aligned}
&\frac{\theta_{1}(2-t)\sqrt{t(1+\delta_{tk})}}{\psi(\theta_{1})}, ~~0<t\leq\frac{2}{3},\\
&\frac{\theta_{2}\sqrt{t(1+\delta_{tk})}}{\varphi(\theta_{2})},  ~~\frac{2}{3}<t\leq1,
\end{aligned}
\right.
~~\widehat{\beta}_{2}=\left\{
\begin{aligned}
&\theta_{1}\sqrt{\frac{(2-t)\delta_{tk}}{\psi(\theta_{1})}}, ~~0<t\leq\frac{2}{3},\\
&\theta_{2}\sqrt{\frac{\delta_{tk}}{\varphi(\theta_{2})}},  ~~\frac{2}{3}<t\leq1,
\end{aligned}
\right.
\end{align*}
with $\theta_{2}=2/t$, and $\psi(\cdot)$ and $\varphi(\cdot)$ being given in \eqref{psi-theta} and \eqref{varphi-theta}, respectively.
\end{theorem}
\begin{remark}
To the best of our knowledge, Theorem \ref{key-Theorem-2} for the first time presents the $tk$-order RIC based coefficient estimate of the RNSP in the case of $0<t\leq1$. This theorem, together with Lemma \ref{Lemma-2}, affirmatively answers under what kind of $tk$-order RIP condition with $t>0$, RNSP will hold. Note that we can also resort to Theorem \ref{key-Theorem-2} to yield a similar result with Theorem \ref{Theorem-1}, and one can find it in Theorem \ref{Theorem-end}. Since the proof Theorem \ref{Theorem-end} is almost same with that of Theorem \ref{Theorem-1}, we here omit it.
\end{remark}
\begin{theorem}\label{Theorem-end}
For any observed vector $\bm{b}=\mathcal{A}(X)+\bm{n}$ with a bounded constraint $\|\bm{n}\|_{2}\leq\epsilon$, if the $tk$-order RIC of $\mathcal{A}$ satisfies \eqref{conditions-step1}, then we have
\begin{align*}
\|X^{\sharp}-X\|_{F}&\leq C_{3}(\widehat{\beta}_{1},\widehat{\beta}_{2})\|X-X_{[k]}\|_{*}+ C_{4}(\widehat{\beta}_{1},\widehat{\beta}_{2}).
\end{align*}
where $X^{\sharp}$ is denoted as the optimal solution of \eqref{NN-Minimization-Unconstrained}.
\end{theorem}
\begin{remark}\label{remark-for-future-work}
Theorem \ref{Theorem-end} and Theorem \ref{Theorem-1} indicate that under a certain $tk$-order RIP condition with $t>0$, the unconstrained RNNM model \eqref{NN-Minimization-Unconstrained} is able to provide a robust recovery performance for any matrix that is not necessary to be exactly low-rank. Note that the previous analysis results on Theorem \ref{Theorem-1} still apply to Theorem \ref{Theorem-end} if one replaces $\beta_{1}$ and $\beta_{2}$ with $\widehat{\beta}_{1}$ and $\widehat{\beta}_{2}$, respectively. On the other hand, one may wonder how the obtained condition \eqref{conditions-step1} performs when compared with the sharp condition \eqref{sharp-deltaTK} established for the constrained NNM model \eqref{NNM-Model}.
\begin{figure}[!htbp]
\begin{center}
\includegraphics[width=0.85\textwidth,height=0.404\textwidth]{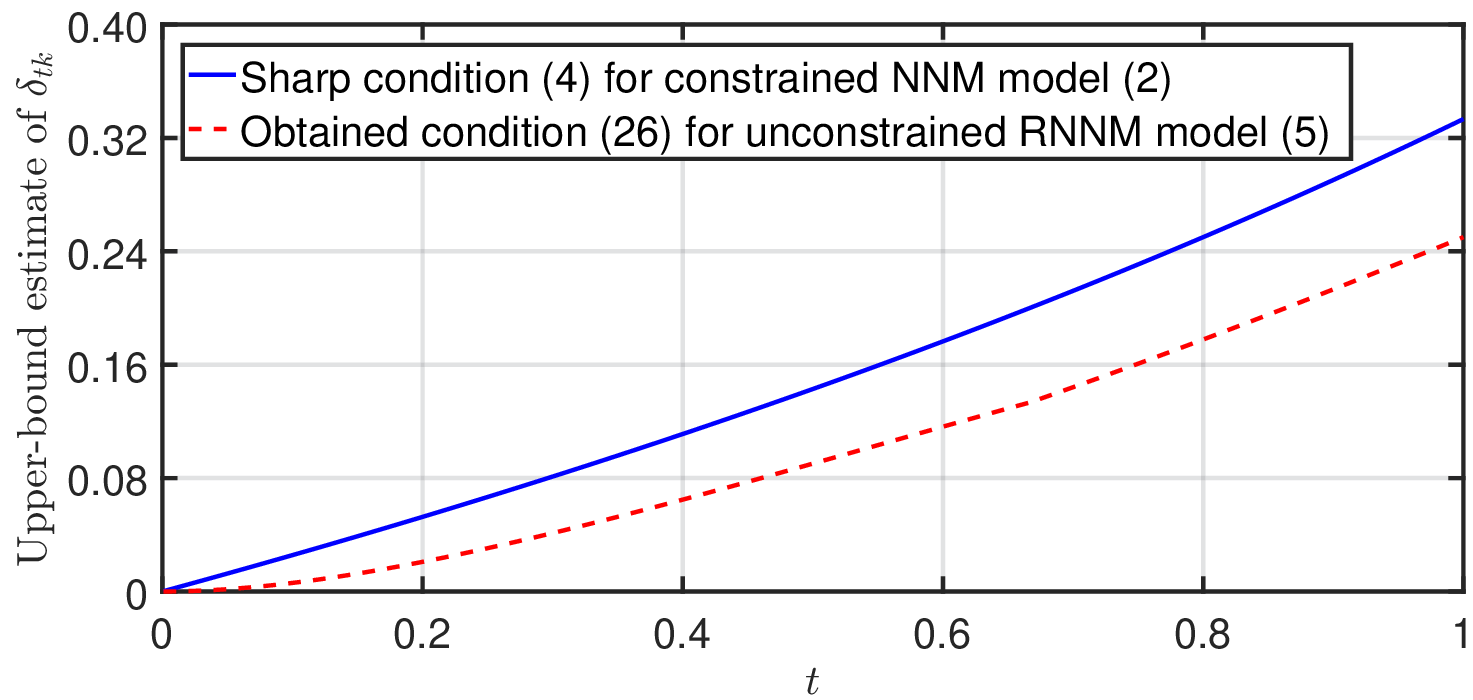}
\caption{\small \textcolor{black}{Comparison between the recovery condition, which takes the form $\delta_{tk}<\delta^{*}$ with $0<t\leq1$, for the constrained model \eqref{NNM-Model} and the unconstrained model \eqref{NN-Minimization-Unconstrained}}}\label{Comparied-conditions}
\end{center}
\end{figure}
Figure \ref{Comparied-conditions} plots the comparison between these two recovery conditions. Unfortunately, our condition \eqref{conditions-step1} is a bit weaker than the sharp condition \eqref{sharp-deltaTK}. Considering the fact that problems of type \eqref{NN-Minimization-Unconstrained} and type \eqref{NNM-Model} are not completely equivalent in both theoretical and applied aspects, it is still an open problem to determine whether the sharp condition \eqref{sharp-deltaTK} for $0<t\leq4/3$ is appropriate for \eqref{NN-Minimization-Unconstrained} or not. According to the established Lemma \ref{Lemma-2} and Theorem \ref{Theorem-1}, it is expected that the condition \eqref{conditions-step1} with $0<t\leq1$ for both Theorem \ref{key-Theorem-2} and Theorem \ref{Theorem-end} can be further improved to the sharp condition \eqref{sharp-deltaTK} with $0<t\leq4/3$. We hope we can solve this problem in the future.
\end{remark}

\begin{proof}[Proof of Theorem \ref{key-Theorem-2}]
\textcolor{black}{
Our proof is inspired by \cite{R-Zhang-A-TIT-2018}. We here only prove the case when $tk$ is an integer since the case when $tk$ is not an integer can be induced easily. Note that for any subset $\Omega\subset[n_{1}]$ with $|\Omega|=k$ and a fixed subset $\widehat{\Omega}=[k]$, $\|H_{\Omega}\|_{F}\leq\|H_{\widehat{\Omega}}\|_{F}$ and $\|H_{\widehat{\Omega}^{c}}\|_{*}\leq\|H_{\Omega^{c}}\|_{*}$ always hold, and hence we will prove
\begin{align}\label{Main-results-2-1END}
\|H_{\widehat{\Omega}}\|_{F}\leq& \widehat{\beta}_{1}\|\mathcal{A}(H)\|_{2}+\widehat{\beta}_{2}\frac{\|H_{\widehat{\Omega}^{c}}\|_{*}}{\sqrt{k}}.
\end{align}
to complete the proof of \eqref{Main-results-2}. We start with denoting the SVD of $H$ as $H=\sum_{i=1}^{n_{1}}\sigma_{i}\bm{a}^{(i)}\left(\bm{c}^{(i)}\right)^{T}$ and $H^{(i)}=\sigma_{i}\bm{a}^{(i)}\left(\bm{c}^{(i)}\right)^{T}$. We denote $a=b=tk/2$ when $tk$ is an even number, and $a=(tk+1)/2$ and $b=(tk-1)/2$ when $tk$ is an odd number. We also denote all the possible index $\Delta_{i}, \Gamma_{j}\subseteq\widehat{\Omega}$ with $|\Delta_{i}|=a$ and $|\Gamma_{j}|=b$, respectively, where $i\in I(|I|=\binom{k}{a})$ and $j\in J(|J|=\binom{k}{b})$. According to Lemma \ref{Key-Combination-Results}, we directly have
\begin{align}\label{Combinication-SUM}
\sum_{i\in I}\|\sigma_{\Delta_{i}}\|_{2}^{2}=\binom{k-1}{a-1}\|\sigma_{\widehat{\Omega}}\|_{2}^{2} \text{~~and~~}\sum_{j\in J}\|\sigma_{\Gamma_{j}}\|_{2}^{2}=\binom{k-1}{b-1}\|\sigma_{\widehat{\Omega}}\|_{2}^{2}.
\end{align}
Besides, we need to introduce the following partition for $\widehat{\Omega}^{c}$, i.e.,
\begin{align}
\Lambda_{3}&=\left\{i\in \widehat{\Omega}^{c}: \sigma_{i}>\frac{\|H_{\widehat{\Omega}^{c}}\|_{*}}{b} \right\},~\Lambda_{4}=\left\{i\in \widehat{\Omega}^{c}: \sigma_{i}\leq\frac{\|H_{\widehat{\Omega}^{c}}\|_{*}}{b} \right\},\label{definition-Lambda12}\\
\Lambda_{5}&=\left\{i\in \widehat{\Omega}^{c}: \sigma_{i}>\frac{\|H_{\widehat{\Omega}^{c}}\|_{*}}{a} \right\}, ~\Lambda_{6}=\left\{i\in \widehat{\Omega}^{c}: \sigma_{i}\leq\frac{\|H_{\widehat{\Omega}^{c}}\|_{*}}{a} \right\}. \label{definition-Lambda34}
\end{align}
By using the similar manipulations as in proof of Lemma \ref{Lemma-2}, we can express $\sigma_{\Lambda_{4}}$ and $\sigma_{\Lambda_{6}}$ as
\begin{align*}
\sigma_{\Lambda_{4}}=\sum_{l}\eta_{l}\bm{u}^{(l)},~\sigma_{\Lambda_{6}}=\sum_{l}\widetilde{\eta}_{l}\bm{v}^{(l)},
\end{align*}
respectively, with $\bm{u}^{(l)}$ and $\bm{v}^{(l)}$ satisfying
\begin{align}
\sum_{l}\eta_{l}\|\bm{u}^{(l)}\|_{2}^{2}\leq (b-|\Lambda_{3}|)\frac{\|H_{\widehat{\Omega}^{c}}\|_{*}^{2}}{b^{2}}\leq \frac{\|H_{\widehat{\Omega}^{c}}\|_{*}^{2}}{b},\label{Z-111-1}\\
\sum_{l}\widetilde{\eta}_{l}\|\bm{v}^{(l)}\|_{2}^{2}\leq (a-|\Lambda_{5}|)\frac{\|H_{\widehat{\Omega}^{c}}\|_{*}^{2}}{a^{2}}\leq \frac{\|H_{\widehat{\Omega}^{c}}\|_{*}^{2}}{a}.\label{Z-111-2}
\end{align}
Let's denote $E^{(l)}=H_{\Lambda_{3}}+ U^{(l)}$ and $ G^{(l)}=H_{\Lambda_{3}}+ U^{(l)}$, where
\begin{align*}
U^{(l)}&=\sum_{i=1}^{n_{1}}\left(\bm{u}^{(l)}\right)_{i}\bm{a}^{(i)}\left(\bm{c}^{(i)}\right)^{T},~
~V^{(l)}=\sum_{i=1}^{n_{1}}\left(\bm{v}^{(l)}\right)_{i}\bm{a}^{(i)}\left(\bm{c}^{(i)}\right)^{T},
\end{align*}
then we have $H_{\Lambda_{4}}=\sum_{l}\eta_{l}U^{(l)}$ and $H_{\Lambda_{6}}=\sum_{l}\widetilde{\eta}_{l}V^{(l)}$. We also need to denote
\begin{align*}
\kappa_{a,b}=&\frac{k-b}{a\binom{k}{a}}\sum_{i\in I}\sum_{l}\eta_{l}\left[a^{2}\left\|\mathcal{A}\left(H_{\Delta_{i}}+\frac{\theta b}{k}E^{(l)}\right)\right\|_{2}^{2}-b^{2}\left\|\mathcal{A}\left(H_{\Delta_{i}}-\frac{\theta a}{k}E^{(l)}\right)\right\|_{2}^{2}\right],\\
\widetilde{\kappa}_{a,b}=&\frac{k-a}{b\binom{k}{b}}\sum_{j\in J}\sum_{l}\widetilde{\eta}_{l}\left[b^{2}\left\|\mathcal{A}\left(H_{\Gamma_{j}}+\frac{\theta a}{k}G^{(l)}\right)\right\|_{2}^{2}-a^{2}\left\|\mathcal{A}\left(H_{\Gamma_{j}}-\frac{\theta b}{k}G^{(l)}\right)\right\|_{2}^{2}\right].
\end{align*}
where $\theta\geq1$ is a number which will be determined later. Since $H_{\Delta_{i}}$, $E^{(l)}$, $H_{\Gamma_{j}}$, $G^{(l)}$ are rank-$a$, -$b$, -$b$, and -$a$, respectively, we can easily induce that any linear combination of $H_{\Delta_{i}}$ and $E^{(l)}$, as well as that of $H_{\Gamma_{j}}$ and $G^{(l)}$, is rank-$tk$. Let's first apply $tk$-order RIP on $\kappa_{a,b}$. This directly gives
\begin{align}\label{kappa-1}
\kappa_{a,b}\geq&\frac{k-b}{a\binom{k}{a}}\sum_{i\in I}\sum_{l}\eta_{l}\left[a^{2}(1-\delta_{tk})\left\|H_{\Delta_{i}}+\frac{\theta b}{k}E^{(l)}\right\|_{2}^{2}-b^{2}(1+\delta_{tk})\left\|H_{\Delta_{i}}-\frac{\theta a}{k}E^{(l)}\right\|_{2}^{2}\right]\nonumber\\
=&\frac{k-b}{a\binom{k}{a}}\sum_{i\in I}\sum_{l}\eta_{l}\bigg\{a^{2}(1-\delta_{tk})\left[\left\|H_{\Delta_{i}}\right\|_{F}^{2}+
\frac{\theta^{2} b^{2}}{k^{2}}\left(\left\|H_{\Lambda_{3}}\right\|_{F}^{2}+ \left\|U^{(l)}\right\|_{F}^{2}\right)\right]\nonumber\\
&-b^{2}(1+\delta_{tk})\left[\left\|H_{\Delta_{i}}\right\|_{F}^{2}+
\frac{\theta^{2} a^{2}}{k^{2}}\left(\left\|H_{\Lambda_{3}}\right\|_{F}^{2}+ \left\|U^{(l)}\right\|_{F}^{2}\right)\right]\bigg\}\nonumber\\
=&\frac{k-b}{k}[(a^{2}-b^{2})-\delta_{tk}(a^{2}+b^{2})]\left\|H_{\widehat{\Omega}}\right\|_{F}^{2}
-\frac{2\theta^{2}(k-a)ab^{2}\delta_{tk}}{k^{2}}\left\|\sigma_{\Lambda_{3}}\right\|_{2}^{2}\nonumber\\
&-\frac{2\theta^{2}(k-b)ab\delta_{tk}}{k^{2}}\|H_{\Omega^{c}}\|_{*}^{2},
\end{align}
where we have used \eqref{Combinication-SUM} and \eqref{Z-111-1} in the last equality. Similarly, we can also apply $tk$-order RIP on $\widetilde{\kappa}_{a,b}$ to get
\begin{align}\label{kappa-2}
\widetilde{\kappa}_{a,b}\geq&\frac{k-a}{k}[(b^{2}-a^{2})-\delta_{tk}(a^{2}+b^{2})]\left\|H_{\widehat{\Omega}}\right\|_{F}^{2}
-\frac{2\theta^{2}(k-b)a^{2}b\delta_{tk}}{k^{2}}\left\|\sigma_{\Lambda_{5}}\right\|_{2}^{2}\nonumber\\
&-\frac{2\theta^{2}(k-a)ab\delta_{tk}}{k^{2}}\|H_{\widehat{\Omega}^{c}}\|_{*}^{2}.
\end{align}
Therefore, combing \eqref{kappa-1} and \eqref{kappa-2} yields
\begin{align}\label{new-add1}
\kappa_{a,b}+\widetilde{\kappa}_{a,b}\geq&[t(a-b)^{2}-(2-t)(a^{2}+b^{2})\delta_{tk}]\left\|H_{\widehat{\Omega}}\right\|_{F}^{2}
-\frac{2\theta^{2}(2-t)ab\delta_{tk}}{k}\|H_{\widehat{\Omega}^{c}}\|_{*}^{2}\nonumber\\
&-\frac{2\theta^{2}ab\delta_{tk}}{k^{2}}\left[b(k-b)\left\|\sigma_{\Lambda_{3}}\right\|_{2}^{2}+a(k-a)\left\|\sigma_{\Lambda_{5}}\right\|_{2}^{2}\right].
\end{align}
Due to the fact that
\begin{align*}
\left\|\sigma_{\Lambda_{3}}\right\|_{2}^{2}\leq\frac{b}{k}\left\|\sigma_{\widehat{\Omega}}\right\|_{2}^{2}=\frac{b}{k}\|H_{\widehat{\Omega}}\|_{F}^{2},~
\left\|\sigma_{\Lambda_{5}}\right\|_{2}^{2}\leq\frac{a}{k}\left\|\sigma_{\widehat{\Omega}}\right\|_{2}^{2}=\frac{a}{k}\|H_{\widehat{\Omega}}\|_{F}^{2},
\end{align*}
we can further induce from \eqref{new-add1} that
\begin{align}\label{24-What}
\kappa_{a,b}+\widetilde{\kappa}_{a,b}\geq&\bigg\{t(a-b)^{2}-(2-t)(a^{2}+b^{2})\delta_{tk}-\frac{2\theta^{2}ab\delta_{tk}}{k^{3}}[b^{2}(k-b)\nonumber\\
&+a^{2}(k-a)]\bigg\}\left\|H_{\widehat{\Omega}}\right\|_{F}^{2}-\frac{2\theta^{2}(2-t)ab\delta_{tk}}{k}\|H_{\widehat{\Omega}^{c}}\|_{*}^{2}\nonumber\\
=&\bigg[t(a-b)^{2}-(2-t)(a^{2}+b^{2})\delta_{tk}-2\theta^{2}(1-t)t^{2}ab\delta_{tk}\nonumber\\
&-\frac{2\theta^{2}(3t-2)a^{2}b^{2}\delta_{tk}}{k^{2}}\bigg]\left\|H_{\widehat{\Omega}}\right\|_{F}^{2}
-\frac{2\theta^{2}(2-t)ab\delta_{tk}}{k}\left\|H_{\widehat{\Omega}^{c}}\right\|_{*}^{2}.
\end{align}
On the other hand, according to the definition of $\kappa_{a,b}$ and $\widetilde{\kappa}_{a,b}$ we can induce that
\begin{align}\label{equlity-kappa1-2}
\kappa_{a,b}+\widetilde{\kappa}_{a,b}=&(a^{2}-b^{2})\underbrace{\left[\frac{k-b}{a\binom{k}{a}}\sum_{i\in I}\|\mathcal{A}\left(H_{\Delta_{i}}\right)\|_{2}^{2}-
\frac{k-a}{b\binom{k}{b}}\sum_{j\in J}\|\mathcal{A}\left(H_{\Gamma_{j}}\right)\|_{2}^{2}\right]}_{\tau}\nonumber\\
&+2\theta abt\left\langle\mathcal{A}\underbrace{\left(\frac{k-b}{a\binom{k}{a}}\sum_{i\in I}H_{\Delta_{i}}+\frac{k-a}{b\binom{k}{b}}\sum_{j\in J}H_{\Gamma_{j}}\right)}_{M}, \mathcal{A}\left(H_{\widehat{\Omega}^{c}}\right) \right\rangle.
\end{align}
To simplify $\tau$ and $M$, from the definition of $\Delta_{i}$ and the SVD of $H$ we can induce that
\begin{align*}
&\sum_{i\in I}\|\mathcal{A}\left(H_{\Delta_{i}}\right)\|_{2}^{2}=\sum_{i\in I}\left\|\sum_{l\in\Delta_{i}}\mathcal{A}\left(H^{(l)}\right)\right\|_{2}^{2}\\
&=\sum_{i\in I}\left[\sum_{l\in\Delta_{i}}\left\|\mathcal{A}\left(H^{(l)}\right)\right\|_{2}^{2}
+\sum_{\substack{p\neq q\\ p,q\in\Delta_{i}}}\left\langle\mathcal{A}\left(H^{(p)}\right),
\mathcal{A}\left(H^{(q)}\right) \right\rangle\right]\\
&=\binom{k-1}{a-1}\sum_{l\in\widehat{\Omega}}\left\|\mathcal{A}\left(H^{(l)}\right)\right\|_{2}^{2}
+\binom{k-2}{a-2}\sum_{\substack{p\neq q\\ p,q\in\widehat{\Omega}}}\left\langle\mathcal{A}\left(H^{(p)}\right),
\mathcal{A}\left(H^{(q)}\right) \right\rangle
\end{align*}
and
\begin{align*}
\sum_{i\in I}H_{\Delta_{i}}=\sum_{i\in I}\sum_{l\in\Delta_{i}}H^{(l)}=\binom{k-1}{a-1}\sum_{l\in \widehat{\Omega}}H^{(l)}=\binom{k-1}{a-1}H_{\widehat{\Omega}},
\end{align*}
where we have used Lemma \eqref{Key-Combination-Results} again. Similarly, we can also get
\begin{align*}
\sum_{j\in J}\|\mathcal{A}\left(H_{\Gamma_{j}}\right)\|_{2}^{2}=&\binom{k-1}{b-1}\sum_{l\in\widehat{\Omega}}\left\|\mathcal{A}\left(H^{(l)}\right)\right\|_{2}^{2}\\
&+\binom{k-2}{b-2}\sum_{\substack{p\neq q\\ p,q\in\widehat{\Omega}}}\left\langle\mathcal{A}\left(H^{(p)}\right)
\mathcal{A}\left(H^{(q)}\right) \right\rangle,
\end{align*}
and $\sum_{j\in J}H_{\Gamma_{j}}=\binom{k-1}{b-1}H_{\widehat{\Omega}}$. Those indicate that
\begin{align*}
\tau=&\frac{a-b}{k}\left[\sum_{l\in\widehat{\Omega}}\left\|\mathcal{A}\left(H^{(l)}\right)\right\|_{2}^{2}+\sum_{\substack{p\neq q\\ p,q\in\widehat{\Omega}}}\left\langle\mathcal{A}\left(H^{(p)}\right)
\mathcal{A}\left(H^{(q)}\right) \right\rangle\right]=\frac{a-b}{k}\left\|\mathcal{A}\left(H_{\widehat{\Omega}}\right)\right\|_{2}^{2}
\end{align*}
and
\begin{align*}
M=\frac{k-b}{a\binom{k}{a}}\sum_{i\in I}H_{\Delta_{i}}+\frac{k-a}{b\binom{k}{b}}\sum_{j\in J}H_{\Gamma_{j}}=(2-t)H_{\widehat{\Omega}},
\end{align*}
and thus we can equivalent write \eqref{equlity-kappa1-2} as
\begin{align}\label{kkk-321}
\kappa_{a,b}+\widetilde{\kappa}_{a,b}=&t(a-b)^{2}\left\|\mathcal{A}\left(H_{\widehat{\Omega}}\right)\right\|_{2}^{2}+2\theta t(2-t) ab\left\langle\mathcal{A}\left(H_{\widehat{\Omega}}\right), \mathcal{A}\left(H_{\widehat{\Omega}^{c}}\right)\right\rangle\nonumber\\
&=t c_{\theta}\left\|\mathcal{A}\left(H_{\widehat{\Omega}}\right)\right\|_{2}^{2}+2\theta t(2-t) ab\left\langle\mathcal{A}\left(H_{\widehat{\Omega}}\right), \mathcal{A}\left(H\right)\right\rangle,
\end{align}
where $c_{\theta}=[(a-b)^{2}-2\theta(2-t) ab]$. Since $\theta\geq1$, $0<t\leq1$, and $a,b$ have been well defined, we can easily know that $c_{\theta}<0$. Note that there are exactly $\binom{k-a}{b}$ sets $\Gamma_{j}$ for a fixed $\Delta_{i}$, exactly $\binom{k-b}{a}$ sets $\Delta_{i}$ for a fixed $\Gamma_{i}$, and exactly $\binom{k-2}{a+b-2}$ sets with $\Delta_{i}\cap\Gamma_{j}=\emptyset$ for two fixed indices $p,q$ with $p\neq q$ and $p,q\in\Delta_{i}\cup\Gamma_{j}$. Therefore, by means of Lemma \ref{Key-Combination-Results} and some simple calculation, we get
\begin{align}\label{KKK-123}
&\frac{1}{\binom{k}{a}\binom{k-a}{b}}\sum_{\Delta_{i}\cap\Gamma_{j}=\emptyset}\left[\left\|\mathcal{A}\left(H_{\Delta_{i}}+H_{\Gamma_{j}}\right)\right\|_{2}^{2}
-\frac{1-t}{ab}\left\|\mathcal{A}\left(bH_{\Delta_{i}}-aH_{\Gamma_{j}}\right)\right\|_{2}^{2}\right]\nonumber\\
&=\frac{1}{\binom{k}{a}\binom{k-a}{b}}\binom{k}{a}\binom{k-a}{b}t^{2}\left\|\mathcal{A}\left(H_{\widehat{\Omega}}\right)\right\|_{2}^{2}
=t^{2}\left\|\mathcal{A}\left(H_{\widehat{\Omega}}\right)\right\|_{2}^{2}.
\end{align}
Moreover, applying $tk$-order RIP on the left-hand side (LHS) of \eqref{KKK-123} and also using \eqref{Combinication-SUM} again, we also have
\begin{align}\label{wa-end}
&\frac{1}{\binom{k}{a}\binom{k-a}{b}}\sum_{\Delta_{i}\cap\Gamma_{j}=\emptyset}\left[\left\|\mathcal{A}\left(H_{\Delta_{i}}+H_{\Gamma_{j}}\right)\right\|_{2}^{2}
-\frac{1-t}{ab}\left\|\mathcal{A}\left(bH_{\Delta_{i}}-aH_{\Gamma_{j}}\right)\right\|_{2}^{2}\right]\geq\nonumber\\
&\frac{1}{\binom{k}{a}\binom{k-a}{b}}\sum_{\Delta_{i}\cap\Gamma_{j}=\emptyset}\left[(1-\delta_{tk})\left\|H_{\Delta_{i}}+H_{\Gamma_{j}}\right\|_{F}^{2}
-\frac{1-t}{ab}(1+\delta_{tk})\left\|bH_{\Delta_{i}}-aH_{\Gamma_{j}}\right\|_{2}^{2}\right]\nonumber\\
&=t[t-(2-t)\delta_{tk}]\|H_{\widehat{\Omega}}\|_{F}^{2}.
\end{align}
Now combining \eqref{kkk-321} and \eqref{KKK-123} yields
\begin{align}\label{Important-1}
&\frac{c_{\theta}}{\binom{k}{a}\binom{k-a}{b}}\sum_{\Delta_{i}\cap\Gamma_{j}=\emptyset}\left\{\left\|\mathcal{A}\left(H_{\Delta_{i}}+H_{\Gamma_{j}}\right)\right\|_{2}^{2}
-\frac{1-t}{ab}\left\|\mathcal{A}\left(bH_{\Delta_{i}}-aH_{\Gamma_{j}}\right)\right\|_{2}^{2}\right\}\nonumber\\
&=t\left(\kappa_{a,b}+\widetilde{\kappa}_{a,b}\right)-2\theta t^{2}(2-t)ab\left\langle\mathcal{A}\left(H_{\widehat{\Omega}}\right), \mathcal{A}\left(H\right)\right\rangle.
\end{align}
As to the LHS of \eqref{Important-1}, we can directly induce from \eqref{wa-end} that
\begin{align}\label{LHS-111}
\text{LHS}\leq c_{\theta}t[t-(2-t)\delta_{tk}]\|H_{\widehat{\Omega}}\|_{F}^{2},
\end{align}
where we have used $c_{\theta}<0$. As to the right-hand side (RHS) of \eqref{Important-1}, we induce from \eqref{24-What} that
\begin{align}\label{RHS-111}
\text{RHS}\geq&t\bigg[t(a-b)^{2}-(2-t)(a^{2}+b^{2})\delta_{tk}-2\theta^{2}(1-t)t^{2}ab\delta_{tk}\nonumber\\
&-\frac{2\theta^{2}(3t-2)a^{2}b^{2}\delta_{tk}}{k^{2}}\bigg]\left\|H_{\widehat{\Omega}}\right\|_{F}^{2}
-\frac{2\theta^{2}t(2-t)ab\delta_{tk}}{k}\left\|H_{\widehat{\Omega}^{c}}\right\|_{*}^{2}\nonumber\\
&-2\theta t^{2}(2-t)ab\sqrt{\frac{1+\delta_{tk}}{t}}\left\|H_{\widehat{\Omega}}\right\|_{F}\left\|\mathcal{A}(H)\right\|_{2}
\end{align}
where, with the aid of \cite[Lemma 4.1]{delta-sharp-Tony-2013}, we have used
\begin{align*}
&\left\langle\mathcal{A}\left(H_{\widehat{\Omega}}\right), \mathcal{A}\left(H\right)\right\rangle\leq\left\|\mathcal{A}\left(H_{\widehat{\Omega}}\right)\|_{2}\|\mathcal{A}\left(H\right)\right\|_{2}
\leq\sqrt{1+\delta_{k}}\left\|H_{\widehat{\Omega}}\|_{F}\|\mathcal{A}\left(H\right)\right\|_{2}\\
&\leq\sqrt{1+\left(\frac{2}{t}-1\right)\delta_{k}}\left\|H_{\widehat{\Omega}}\|_{F}\|\mathcal{A}\left(H\right)\right\|_{2}
\leq\sqrt{\frac{1+\delta_{tk}}{t}}\left\|H_{\widehat{\Omega}}\|_{F}\|\mathcal{A}\left(H\right)\right\|_{2}.
\end{align*}
Therefore, we can induce from \eqref{LHS-111} and \eqref{RHS-111} that
\begin{align*}
&\bigg\{t(a-b)^{2}-(2-t)(a^{2}+b^{2})\delta_{tk}-c_{\theta}[t-(2-t)\delta_{tk}]-2\theta^{2}(1-t)t^{2}ab\delta_{tk}\\
&-\frac{2\theta^{2}(3t-2)a^{2}b^{2}}{k^{2}}\delta_{tk}\bigg\}\left\|H_{\widehat{\Omega}}\right\|_{F}^{2}-2\theta t(2-t)ab\sqrt{\frac{1+\delta_{tk}}{t}}\left\|\mathcal{A}(H)\right\|_{2}\left\|H_{\widehat{\Omega}}\right\|_{F}\\
&-\frac{2\theta^{2}(2-t)ab\delta_{tk}}{k}\left\|H_{\widehat{\Omega}^{c}}\right\|_{*}^{2}\leq0,
\end{align*}
which is also equivalent to
\begin{align}\label{Final-inequlity}
&f(\theta)\left\|H_{\widehat{\Omega}}\right\|_{F}^{2}-\theta (2-t)\sqrt{t(1+\delta_{tk})}\left\|\mathcal{A}(H)\right\|_{2}\left\|H_{\widehat{\Omega}}\right\|_{F}
-\frac{\theta^{2}(2-t)\delta_{tk}}{k}\left\|H_{\widehat{\Omega}^{c}}\right\|_{*}^{2}\leq0,
\end{align}
where $f(\theta)=(2-t)\{t\theta-[1+(2-t)\theta]\delta_{tk}\}-\theta^{2}(1-t)t^{2}\delta_{tk}-\theta^{2}(3t-2)ab\delta_{tk}/k^{2}$.
}

\textcolor{black}{
\textbf{Case 1}: $0<t\leq2/3$. In this case, we know from $f(\theta)$ that
\begin{align}\label{psi-theta}
f(\theta)\geq&(2-t)\{t\theta-[1+(2-t)\theta]\delta_{tk}\}-\theta^{2}(1-t)t^{2}\delta_{tk}\nonumber\\
=&(2-t)t\theta-\left\{(2-t)[1+(2-t)\theta]+(1-t)t^{2}\theta^{2}\right\}\delta_{tk}\triangleq\psi(\theta).
\end{align}
To make sure $f(\theta)>0$, we only need to set $\psi(\theta)>0$, i.e.,
\begin{align*}
\delta_{tk}<\frac{(2-t)t\theta}{(2-t)[1+(2-t)\theta]+(1-t)t^{2}\theta^{2}}=\frac{(2-t)t}{(2-t)^{2}+(2-t)/\theta+(1-t)t^{2}\theta}.
\end{align*}
To obtain as large an upper bound as possible, we need to set $\theta=\theta_{1}=\frac{1}{t}\sqrt{\frac{2-t}{1-t}}$, and thus the largest upper bound with respect to $\theta$ will take the form
\begin{align}\label{neewly-condition-11}
\delta_{tk}<\frac{(2-t)t}{(2-t)^{2}+2t\sqrt{(2-t)(1-t)}}.
\end{align}
One can easily check that $\theta_{1}\geq1$ holds for any $0<t\leq2/3$. Based on the above settings, we can further know from \eqref{Final-inequlity} that
\begin{align*}
&\psi(\theta_{1})\left\|H_{\widehat{\Omega}}\right\|_{F}^{2}-\theta_{1} (2-t)\sqrt{t(1+\delta_{tk})}\left\|\mathcal{A}(H)\right\|_{2}\left\|H_{\widehat{\Omega}}\right\|_{F}
-\frac{\left(\theta_{1}\right)^{2}(2-t)\delta_{tk}}{k}\left\|H_{\widehat{\Omega}^{c}}\right\|_{*}\leq0,
\end{align*}
and hence get
\begin{align*}
\left\|H_{\widehat{\Omega}}\right\|_{F}\leq\frac{\theta_{1}(2-t)\sqrt{t(1+\delta_{tk})}}{\psi(\theta_{1})}\left\|\mathcal{A}(H)\right\|_{2}
+\theta_{1}\sqrt{\frac{(2-t)\delta_{tk}}{\psi(\theta_{1})}}\frac{\left\|H_{\widehat{\Omega}^{c}}\right\|_{*}}{\sqrt{k}},
\end{align*}
which is the desired \eqref{Main-results-2} for $0<t\leq2/3$. To guarantee that \eqref{Main-results-2} satisfies the RNSP with $\widehat{\beta}_{1}>0$ and $0<\widehat{\beta}_{2}<1$, we have to set $\widehat{\beta}_{2}=\theta_{1}\sqrt{(2-t)\delta_{tk}/\psi(\theta_{1})}<1$, i.e.,
\begin{align}\label{neewly-condition-12}
\delta_{tk}<\frac{(2-t)t}{(2-t)^{2}+2t\sqrt{(2-t)(1-t)}+(2-t)\theta_{1}}.
\end{align}
Obviously, the obtained condition \eqref{neewly-condition-12} is the desired condition \eqref{conditions-step1} for $0<t\leq2/3$, which is also included in \eqref{neewly-condition-11}.
}

\textcolor{black}{
\textbf{Case 2}: $2/3<t\leq1$. In this case, by using $ab\leq t^{2}k^{2}/4$, we can also know from $f(\theta)$ that
\begin{align}\label{varphi-theta}
f(\theta)\geq&(2-t)\{t\theta-[1+(2-t)\theta]\delta_{tk}\}-\theta^{2}(1-t)t^{2}\delta_{tk}-\frac{\theta^{2}t^{2}(3t-2)}{4}\delta_{tk},\nonumber\\
=&(2-t)\underbrace{\left\{t\theta-\left[1+(2-t)\theta+t^{2}\theta^{2}/4\right]\delta_{tk}\right\}}_{\triangleq\varphi(\theta)}.
\end{align}
To make sure that $f(\theta)>0$ and $\delta_{tk}$ has as large an upper bound as possible, By using the similar manipulations as in Case 1, we select $\theta=\theta_{2}=2/t$, and thus get $\delta_{tk}<t/2$. Obviously, $\theta_{2}\geq1$ holds for any $2/3<t\leq1$. Furthermore, we can also know from \eqref{Final-inequlity} that
\begin{align*}
&\varphi(\theta_{2})\left\|H_{\widehat{\Omega}}\right\|_{F}^{2}-\theta_{2} \sqrt{t(1+\delta_{tk})}\left\|\mathcal{A}(H)\right\|_{2}\left\|H_{\widehat{\Omega}}\right\|_{F}
-\frac{\left(\theta_{2}\right)^{2}\delta_{tk}}{k}\left\|H_{\widehat{\Omega}^{c}}\right\|_{*}\leq0,
\end{align*}
and hence get
\begin{align*}
\left\|H_{\widehat{\Omega}}\right\|_{F}\leq\frac{\theta_{2}\sqrt{t(1+\delta_{tk})}}{\varphi(\theta_{2})}\left\|\mathcal{A}(H)\right\|_{2}
+\theta_{2}\sqrt{\frac{\delta_{tk}}{\varphi(\theta_{2})}}\frac{\left\|H_{\widehat{\Omega}^{c}}\right\|_{*}}{\sqrt{k}},
\end{align*}
which is the desired \eqref{Main-results-2} for $2/3<t\leq1$. Similarly, to enforce \eqref{Main-results-2} to obey the RNSP with $\widehat{\beta}_{1}>0$ and $0<\widehat{\beta}_{2}<1$, we have to set $\widehat{\beta}_{2}=\theta_{2}\sqrt{\delta_{tk}/\varphi(\theta_{2})}<1$, i.e.,
\begin{align*}
\delta_{tk}<\frac{t^{2}}{2(t+1)},
\end{align*}
which is the desired condition \eqref{conditions-step1} for $2/3<t\leq1$. Combining Case 1 and Case 2, we obtain the desired \eqref{Main-results-2-1END}, and thus establish the results showed in Theorem \ref{key-Theorem-2}.
}
\end{proof}

\section{Conclusion and future work}\label{SECTION-5}
This paper has considered the robust matrix recovery from the unconstrained RNMM model. First, equipped with the powerful $tk$-order RIP tool for $t>0$, we developed a family of $tk$-order RIC based coefficient estimates for the RNSP. To the best of our knowledge, the obtained RNSP results in the case of $0<t\leq1$ have not been explored before. Furthermore, by mean of these RNSP results, some upper-bound estimates of error were established for the unconstrained RNMM model to guarantee the robust matrix recovery. As we have pointed out in Remark \ref{remark-for-future-work}, one of our future work will focus on extending the condition \eqref{conditions-step1} with $0<t\leq1$ to the sharp condition \eqref{sharp-deltaTK} with $0<t\leq4/3$. Besides, determining a proper from the theoretical aspect for the unconstrained RNMM model will be another future work.

\section*{Acknowledgement}
The authors would like to thank the editors and referees for their valuable comments that improve the presentation of this paper.

\end{document}